\documentclass[10pt]{article}
\usepackage{amsfonts}
\usepackage{amsmath}
\usepackage{mathrsfs}
\usepackage{mathrsfs,amscd,amssymb,amsthm,amsmath,bm,graphicx,psfrag,subfigure}

\makeatletter

\renewcommand{\@seccntformat}[1]{{\csname the#1\endcsname}{\normalsize .}\hspace{.5em}}
\makeatother

\def \[{\begin{equation}}
\def \]{\end{equation}}

\def \dist{{\rm dist}}
\newtheorem{thm}{Theorem}[section]

\newtheorem{claim}{Claim}
\newtheorem{lem}[thm]{Lemma}

\newenvironment{wst}
{\setlength{\leftmargini}{1.5\parindent}
 \begin{itemize}
 \setlength{\itemsep}{-1.1mm}}
{\end{itemize}}

\begin{document}
\setlength{\baselineskip}{15pt}
\begin{center}{\Large \bf On the spectral moment of graphs with $k$ cut edges\footnote{Financially supported by the National Natural
Science Foundation of China (Grant Nos. 11071096, 11271149) and the Special Fund for Basic Scientific Research of Central Colleges (CCNU11A02015).}}

\vspace{4mm}

{\large Shuchao Li$^{a,}$\footnote{E-mail: lscmath@mail.ccnu.edu.cn (S.C.
Li), 425333559@qq.com (H. Zhang)},\ \ Huihui Zhang$^a$,\ \ Minjie Zhang$^b$}\vspace{2mm}

$^a$Faculty of Mathematics and Statistics,  Central China Normal
University, Wuhan 430079, P.R. China\vspace{1mm}

$^b$School of Mathematics and Physics, Hubei Institute of Technology, Huangshi 435003, P.R. China
\end{center}
\noindent {\bf Abstract}: Let $A(G)$ be the adjacency matrix of a graph $G$ with $\lambda_{1}(G)$, $\lambda_{2}(G)$, $\dots$, $\lambda_{n}(G)$ being its eigenvalues in non-increasing order. Call the number $S_k(G):=\sum_{i=1}^{n}\lambda_{i}^k(G)\, (k=0,1,\dots,n-1)$ the $k$th spectral moment of $G$. Let $S(G)=(S_0(G),S_1(G),\dots,S_{n-1}(G))$ be the sequence of spectral moments of $G$. For two graphs $G_1$ and $G_2$, we have $G_1\prec_sG_2$ if $S_i(G_1)=S_i(G_2)\, (i=0,1,\dots,k-1)$ and $S_k(G_1)<S_k(G_2)$ for some $k\in \{1,2,\emph{}\dots,n-1\}$. Denote by $\mathscr{G}_n^k$ the set of connected $n$-vertex graphs with $k$ cut edges. In this paper, we determine the first, the second, the last and the second
last graphs, in an $S$-order, among $\mathscr{G}_n^k$, respectively.

\vspace{2mm} \noindent{\it Keywords}: Spectral moment; Cut edge; Clique

\vspace{2mm}


 {\setcounter{section}{0}
\section{\normalsize Introduction}\setcounter{equation}{0}
All graphs considered here are finite, simple and connected. Undefined terminology and notation may be referred to \cite{D-I}.
Let $G=(V_G,E_G)$ be a simple undirected graph with $n$ vertices. $G-v$, $G-uv$ denote the graph obtained from $G$ by deleting vertex $v \in V_G$, or edge
$uv \in E_G$, respectively (this notation is naturally extended if more than one vertex or edge is deleted). Similarly,
$G+uv$ is obtained from $G$ by adding an edge $uv \not\in E_G$. For $v\in V_G$, let
$N_G(v)$ (or $N(v)$ for short) denote the set of all the adjacent vertices of $v$ in $G$ and $d_G(v)=|N_G(v)|$, and $\dist_G(u,v)$ is the distance between $u$ and $v$. For an edge subset $E'$ of $G$, denoted by $G[E']$ the subgraph induced by $E'$.
A cut edge in a connected graph $G$ is an edge whose deletion breaks the graph into two components. Let $\mathscr{G}_n^k$ be the set of all $n$-vertex graphs, each of which contains $k$ cut edges.

Let $A(G)$ be the adjacency matrix of a graph $G$ with $\lambda_1(G),\lambda_2(G),\dots,\lambda_n(G)$ being its eigenvalues in non-increasing order. The number $\sum_{i=1}^n\lambda_i^k(G)\, (k=0,1,\dots,n-1)$ is called the $k$th spectral moment of $G$, denoted by $S_k(G)$. Let $S(G)=(S_0(G), S_1(G), \dots, S_{n-1}(G))$ be the sequence of spectral moments of $G$. For two graphs $G_1, G_2$, we shall write $G_1=_sG_2$ if $S_i(G_1)=S_i(G_2)$ for $i=0,1,\dots,n-1$. Similarly, we have $G_1\prec_sG_2\, (G_1$ comes before $G_2$  in an $S$-order) if for some $k\, (1\leq k\leq {n-1})$, we have  $S_i(G_1)=S_i(G_2)\, (i=0,1,\dots,k-1)$  and $S_k(G_1)<S_k(G_2)$. We shall also write $G_1\preceq_sG_2$ if $G_1\prec_sG_2$ or $G_1=_sG_2$.  $S$-order has been used in producing graph catalogs (see \cite{C-G}), and for a more general setting of spectral moments one may be referred to \cite{C-R-S1}.

Investigation on $S$-order of graphs attracts more and more researchers' attention. Cvetkovi\'c and Rowlinson \cite{C-R-S3} studied
the $S$-order of trees and unicyclic graphs and characterized the first and the last graphs, in an $S$-order,
of all trees and all unicyclic graph with given girth, respectively. Chen, Liu and Liu \cite{C-L-L} studied the lexicographic ordering
by spectral moments ($S$-order) of unicyclic graph with a given girth. Wu and Fan \cite{D-F} determined the first
and the last graphs, in an $S$-order, of all unicyclic graphs and bicyclic graphs, respectively. Pan et al. \cite{X-F-P}
gave the first $\sum_{k=1}^{\lfloor\frac{n-1}{3}\rfloor}(\lfloor\frac{n-k-1}{2}\rfloor-k+1)$ graphs apart from an $n$-vertex path, in an $S$-order, of all trees with $n$ vertices. Wu and Liu \cite{D-M} determined the last $\lfloor\frac{d}{2}\rfloor+1$, in an $S$-order,
among all $n$-vertex trees of diameter $d\, (4 \le d \le n-3)$. Pan et al. \cite{B-Z1} identified the
last and the second last graphs, in an $S$-order, of quasi-trees. Hu, Li and Zhang \cite{Li-H} studied the spectral moments of graphs with given clique number and chromatic number, respectively. Li and Song \cite{Li-S} identified the last $n$-vertex tree with a given degree sequence in an $S$-order. Consequently, the last trees in an $S$-order in the sets of all trees of order $n$ with the largest degree, the leaves number, the independence number and the matching number was also determined, respectively.

In light of the information available from the related results on the spectral moments of graphs, it is natural to consider this problem on some other class of graphs, and the connected graphs with $k$ cut edges are a reasonable starting point for such a investigation. The $n$-vertex connected graphs with $k$ cut edges have been considered in different fields  \cite{B-Z10,2,3,4}, whereas to our best knowledge, the spectral moments of graphs in $\mathscr{G}_n^k$ were, so far, not considered. Here, we identified the first, the second, the last and the second
last graphs, in an $S$-order, among $\mathscr{G}_n^k$, respectively.

Throughout the text we denote by $P_n, K_{1,n-1}, C_n$ and $K_n$ the path, star, cycle and complete graph on $n$ vertices, respectively. Let $K_{1,n-1}^*$ be a graph obtained from a star $K_{1,n-1}$ by attaching a leaf to one leaf of $K_{1,n-1}$, $U_n$ be a graph obtained from $C_{n-1}$ by attaching a leaf to one vertex of $C_{n-1}$, and $B_4, B_5$ be two graphs  obtained from two cycle $C_3, C_3'$ of length 3 by identifying one edge of  $C_3$ with one edge of $C_3'$ and identifying one vertex of $C_3$ with one vertex  of $C_3'$, respectively; see Fig. 1.
\begin{figure}[h!]
\begin{center}
\psfrag{a}{$u_1$}
\psfrag{b}{$u_2$}\psfrag{e}{$K_{a_1}$}
\psfrag{f}{$K_{a_2}$}
\psfrag{g}{$K_{a_{k-1}}$}
\psfrag{h}{$K_{a_k}$}
\psfrag{c}{$u_{k-1}$}
\psfrag{d}{$u_k$}
\psfrag{x}{$u_0$}
\psfrag{y}{$v_0$}
\psfrag{A}{$K_{a_0}$}
\psfrag{9}{$K(a_0,\{a_1, a_2, \dots, a_k\})$}
\psfrag{1}{$B_4$}
\psfrag{2}{$B_5$}
\psfrag{3}{$K_6^1$}
\psfrag{4}{$K_6^2$}
\psfrag{5}{$K_6^3$}
\psfrag{6}{$P_6^1$}
\psfrag{7}{$P_6^2$}
\psfrag{8}{$P_6^3$}
\includegraphics[width=120mm]{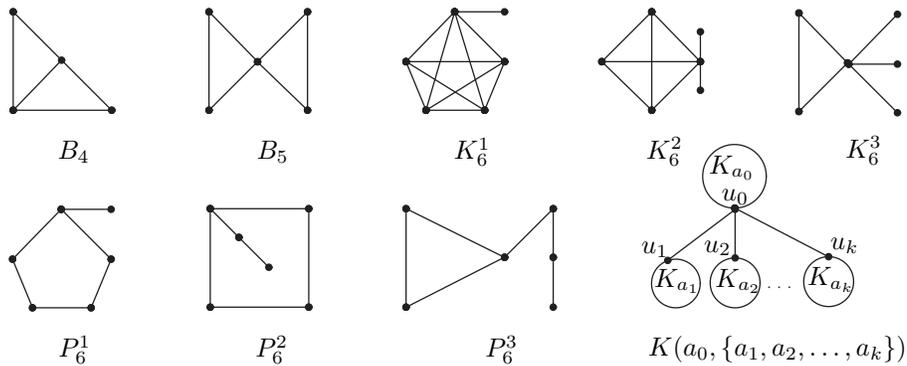}\\
\caption{Graphs $B_4, B_5, K_6^1, K_6^2, K_6^3, P_6^1, P_6^2, P_6^3$ and $K(a_0,\{a_1, a_2, \dots, a_k\}).$}
\end{center}
\end{figure}

The graph $K_n^k$ is an $n$-vertex graph obtained by attaching $k$ pendant vertices to one vertex of $K_{n-k}$. The graph $P_n^k$ is a graph obtained by identifying one end-vertex of $P_{k+1}$ with one vertex of $C_{n-k}$. For example, for $n=6, K_6^0=K_6, K_6^5  $ is a star, $P_6^0=C_6$ and $K_6^1, K_6^2, K_6^3, P_6^1,  P_6^2, P_6^3$ are depicted in Fig. 1. In general, $K_n^0=K_n, K_n^{n-1}$ is star $K_{1,n-1}$, $K_n^{n-2}\cong K_n^{n-1}$ and $P_n^0=C_n$. Let  $K(a_0,\{a_1,a_2,\dots,a_k\})$ be a graph obtained from $K_{1,k}$ by replacing each $u_i\in V_{K_{1,k}}$ by a clique $K_{a_i}\, (a_i\geq1, i=0,1,2,\dots,k)$; see Fig. 1. Denote
$$
 \mathscr{K}_n^k=\left\{K(a_0,\{a_1,a_2,\dots,a_k\}):a_i\geq1(0\leq i\leq k),\,\sum_{i=0}^ka_i=n\right\}.
$$
Let $F$ be a graph. An $F$-subgraph of $G$ is a subgraph of $G$ which is isomorphic to the graph $F$. Let $\phi_{G}(F)$ (or $\phi(F)$) be the number of all $F$-subgraph of $G$.
\begin{lem} {\rm(see \cite{C-R-S2})}
The $k$th spectral moment of $G$ is equal to the number of closed walks of length $k$.
\end{lem}
\begin{lem}
For every graph $G$, we have
\begin{wst}
\item[{\rm (i)}] $S_4(G)=2\phi(P_2)+4\phi(P_3)+8\phi(C_4)$ {\rm(see \cite{C-G});}
\item[{\rm (ii)}] $S_5(G)=30\phi(C_3)+10\phi(U_4)+10\phi(C_5)$ {\rm(see\cite{D-M});}
\item[{\rm (iii)}]
     $S_6(G)=2\phi(P_2)+12\phi(P_3)+6\phi(P_4)+12\phi(K_{1,3})+12\phi(U_5)+36\phi(B_4)+24\phi(B_5)+24\phi(C_3)+48\phi(C_4)+12\phi(C_6)
\, ${\rm(see\cite {C-L-L}).}
\end{wst}
\end{lem}
\begin{lem}[\cite{C-R-S2}]
Given a connected graph $G$, $S_0(G)=n, S_1(G)=l, S_2(G)=2m, S_3(G)=6t$, where $n, l, m, t$ denote the number of vertices, the number of loops, the number of edges and the number of triangles contained in $G$, respectively.
\end{lem}
\begin{lem}[\cite{C-R-S3}]\label{lem4}
In an $S$-order of the $n$-vertex unicyclic graphs with girth $g$, the first graph is $U_n^g$ which is obtained by the coalescence of a cycle $C_g$  with
a path $P_{n-g+1}$ at one of its end-vertices.
\end{lem}
\section{\normalsize The last and the second last graphs in an $S$-order among $\mathscr{G}_n^k$}\setcounter{equation}{0}
In this section, we will determine the last two graphs, in an $S$-order, among $\mathscr{G}_n^k$.
Let $\Bbb E=\{e_1,e_2,\dots,e_k\}$ be the set of the cut edges of $G\in\mathscr{G}_n^k$. Note that $S_2(G)=2|E_G|$, hence $S_2(G+e)>S_2(G)$. By Lemma 1.3, in order to determine the last graph in an $S$-order among $\mathscr{G}_n^k$, it suffices to choose graph $G\in\mathscr{G}_n^k\emph{}$ such that its $S_2(G)$ is as large as possible. So we can have the following assumption throughout this section.

\vspace{2mm}\noindent {\bf Assumption 0.}  Each component of $G-\Bbb E$ is a clique.

\begin{thm}
Of all the connected graphs with $n$ vertices and $k$ cut edges, the last graph in an $S$-order is obtained uniquely at $K_n^k$.
\end{thm}

\begin{proof}

If $k=0$, then by Assumption 0 we have $\mathscr{G}_n^0=\{K_n\}$, our result holds immediately. Therefore we may assume that $k\geq1$. Again by Assumption 0, we can denote the components of $G-\Bbb E$ by $K_{a_0},K_{a_1},\dots,K_{a_k},\, a_0+a_1+\dots+a_k=n$. Assume, without loss of generality, that $a_0\geq a_1\geq a_2\geq\dots\geq a_k\geq1$.

Let $V_i=\{v\in V_{K_{a_i}}\!:\, v$ is an end-vertex of a cut edge of $G$\}. Choose $G\in\mathscr{G}_n^k$  such that $G$ is as large as possible under the order $\preceq_s$. In order to complete the proof, it suffices to show the following facts.

\vspace{2mm}\noindent {\bf Fact 1}.  $|V_i|=1$ for $i=0,1,2,\ldots, k.$
\begin{proof}
Suppose to the contrary that there exists $i\in \{0,1,2,\ldots, k\}$ such that $|V_i|>1$. Let $u,u' \in V_{a_i}$, both $u$ and $u'$ are end-vertices of the cut edges of $G$. Denote  $N_G(u)\backslash N_{K_{a_i}}(u)=\{w_1,w_2,\dots,w_s\}$ and $N_G(u')\backslash N_{K_{a_i}}(u')=\{z_1,z_2,\dots,z_l\}$. It is routine to check that  $s\geq1,\, l\geq 1$. Let
  $$G^*=G-\{u'z_1,u'z_2,\dots,u'z_l\}+\{uz_1,uz_2,\dots,uz_l\},
  $$
  then  $G^* \in \mathscr{G}_n^k$.

On the one hand, $S_i(G)=S_i(G^*)$ for $i=0, 1, 2, 3.$ On the other hand, $\phi_G(P_2)=\phi_{G^*}(P_2),  \phi_G(C_4)=\phi_{G^*}(C_4)$, hence by Lemma 1.2(i),
\begin{eqnarray*}
 S_4(G)-S_4(G^*)  = 4(\phi_G(P_3)-\phi_{G^*}(P_3))=4\left({s\choose 2}+{l\choose 2}-{s+l\choose 2}\right)=-4sl<0.
  \end{eqnarray*}
which implies that $G\prec_sG^*$, a contradiction. Therefore $|V_i|=1$ for $0\leq i \leq k $.
\end{proof}

By Fact 1, we can assume that $V_i=\{u_i\}$ for $i=0,1,2,\ldots, k$.

\vspace{2mm}\noindent {\bf Fact 2}.  $G\in \mathscr{K}_n^k$.\vspace{2mm}

\begin{proof}
If not, then there exists a cut edge $u_0u_i\in \Bbb E$ such that $u_i$ is an end-vertex of another cut edge(s). Let
$$
|N_G(u_i)\setminus (N_{K_{a_i}}(u_i)\cup \{u_0\})|=l, \ \ \ |N_G(u_0)\setminus (N_{K_{a_0}}(u_0)\cup \{u_i\})|=s.
$$
It is straightforward to check that $l\ge 1$ and $s\ge 0.$

First consider that $s\geq1$. In this case, let
$$
    G^*=G-\{u_iz: \, z\in N_G(u_i)\setminus (N_{K_{a_i}}(u_i)\cup \{u_0\})\}+\{u_0z: \, z\in N_G(u_i)\setminus (N_{K_{a_i}}(u_i)\cup \{u_0\})\}.
$$
It is easy to see that $G^*\in \mathscr{G}_n^k$. Note that $S_i(G)=S_i(G^*)$ for $i=0,1,2,3$ and $\phi_G(P_2)=\phi_{G^*}(P_2)$, $\phi_G(C_4)=\phi_{G^*}(C_4)$, hence
    $$
    S_4(G)-S_4(G^*)  = 4(\phi_G(P_3)-\phi_{G^*}(P_3))=4(l(a_i-1)-l(a_0-1)-ls)=4l(a_i-a_0-s)<0,
    $$
which implies that $G\prec_sG^*$, a contradiction.

Now consider that $s=0$. In this case, there exists a cut edge $u_iu_j\in \Bbb E$ such that $u_j$ is an end-vertex of another cut edge(s).
Let
$
|N_G(u_j)\setminus (N_{K_{a_j}}(u_j)\cup \{u_i\})|=p.
$
It is straightforward to check that $p\ge 1$. Let
\begin{eqnarray*}
    G^*&=&G-\{u_jz: \, z\in N_G(u_j)\setminus (N_{K_{a_j}}(u_j)\cup \{u_i\})\}+\{u_0z: \, z\in N_G(u_j)\setminus (N_{K_{a_j}}(u_j)\cup \{u_i\})\}\\
       && -\{u_iw: \, w\in N_G(u_i)\setminus (N_{K_{a_i}}(u_i)\cup \{u_0\})\}+\{u_0w: \, w\in N_G(u_i)\setminus (N_{K_{a_i}}(u_i)\cup \{u_0\})\}.
\end{eqnarray*}
It is easy to see that $G^*\in \mathscr{G}_n^k$. Note that $S_i(G)=S_i(G^*)$, $i=0,1,2,3$ and $\phi_G(P_2)=\phi_{G^*}(P_2)$, $\phi_G(C_4)=\phi_{G^*}(C_4)$. Hence,
\begin{eqnarray*}
S_4(G)-S_4(G^*)  &=&  4(\phi_G(P_3)-\phi_{G^*}(P_3))\\
                 &=& 4l(a_i-1)+p(a_j-1)-4p(l-1)-4p-4(l+p)(a_0-1)\\
                 &=&4[l(a_i-a_0)+p(a_j-a_0)-pl]<0.
\end{eqnarray*}
The last inequality follows from $a_i\le a_0, a_j\le a_0$ and $pl>0.$ Hence, we obtain that $G\prec_sG^*$, a contradiction. Therefore $G\in \mathscr{K}_n^k$.
\end{proof}

By Fact 2, we can assume that $u_0u_j\in \Bbb E$, $1\leq j \leq k$.

\vspace{2mm}\noindent {\bf Fact 3}. $a_1=a_2=\cdots=a_k=1$.
\begin{proof}
Assume to the contrary that there exists a $j\in \{1,2,\ldots,k\}$ such that $a_j>1$. By Fact 2, we have $G=K(a_0,\{a_1,\ldots,a_{j-1},a_j, a_{j+1},\ldots,a_k\})$. Now we consider $G^*=K(a_0+a_j-1,\{a_1,\ldots,a_{j-1},1, a_{j+1},\linebreak \ldots,a_k\})$. It is easy to see that $G^*\in \mathscr{K}_n^k.$

Note that $S_i(G)=S_i(G^*),\,i=0,1$ and
$$
  S_2(G)-S_2(G^*)=2(a_j-1)-2(a_j-1)a_0=2(a_j-1)(1-a_0)<0,
$$
i.e., $ G\prec_s G^*$, a contradiction. Therefore $a_j=1$ for $j=1,2,\ldots, k$.
 \end{proof}
In the view of Fact 3, we have $a_0=n-k$. Hence, $G=K(n-k,\{1,1,\ldots,1\})$, i.e., $G\cong K_n^k$, as desired.
\end{proof}


In the rest of this section, we are to determine last graph in an $S$-order among $\mathscr{G}_n^k\backslash K_n^k$.
Delete an edge, say $xy$, from $K_n$ and denote the resultant graph by $G_1$. Let $G_2$ be a graph obtained from $G_1$ by attaching a pendant vertex to one vertex, say $r$, of $G_1$ with $r\not=x, y$. Let $G_3=K(n-k,\{\underbrace{1,1,\dots,1}_k\})-uw+vw$, where $uw$ is a cut edge and $u, v$ are two different vertices in $V_{K_{n-k}}$.

Based on Lemma 1.3, it is easy to see that among $\mathscr{G}_n^0,$ $K_n$ (resp. $G_1$) is the last (resp. the second last) graph in an $S$-order, while among $\mathscr{G}_n^1$ with $n\ge 5$, based on $S_2(G)$, the second last graph in an $S$-order must be a graph obtained from $K_n^1$ by deleting a non-cut edge, say $e$, from $K_n^1$. Denote the resultant graph by $G'$ if $e$ has a common vertex with the cut edge in $K_n^1$ and by $G_2$ otherwise. Note that $S_i(G_2)=S_i(G')$ for $i=0,1,2,3$ and
$\phi_{G_2}(P_2)=\phi_{G'}(P_2),\, \phi_{G_2}(C_4)=\phi_{G'}(C_4),$ hence by Lemma 1.2(i)
$$
S_4(G')-S_4(G_2)=4(\phi_{G'}(P_3)-\phi_{G_2}(P_3))=-4<0,
$$
i.e., $G'\prec_sG_2$, Hence, among $\mathscr{G}_n^1$ with $n\ge 5$, $G_2$ is the second last graph in an $S$-order. In what follows we only consider $k\geq2$.

\begin{thm}
Among $\mathscr{G}_n^k$ with $2\le k\le n-1$, the second last graph in an $S$-order is obtained uniquely at $G_3$ if $k\in \{2,3,\ldots, n-2\}$ and at $K_{1,n-1}^*$ otherwise, where $G_3$ is defined as above.
\end{thm}
\begin{proof}
Choose $G\in \mathscr{G}_n^k\setminus\{K_n^k\}$ such that it is as large as possible according to $\preceq_s$.
Denote the components of $G-\Bbb E$ by $U_0, U_1, U_2, \ldots, U_k$. We are to show that each of the components is a complete graph. In fact, if there exists a $U_i$ which is not a complete graph, i.e., $U_i$ contains two vertices $x, y$ satisfying $xy\not\in E_{U_i}$. Let $G'=G+xy$. If $G'\not\cong K_n^k$, it is easy to see that $G\prec_s G'$, a contradiction. If $G'\cong K_n^k$, then either $x$ or $y$ is not an end-vertex of a cut edge of $G$. Without loss of generality, assume that  $x$ is not an end-vertex of a cut edge of $G$, delete a cut edge of $G'$ and connect the isolated vertex with $x$ by an edge; denote the resultant graph by $G''$. Then we have $S_0(G)=S_0(G''), S_1(G)=S_1(G'')$ and $S_2(G)<S_2(G'').$ Hence, $G\prec_s G''$, a contradiction. Therefore, we may denote the components of $G-\Bbb E$ by $K_{a_0},K_{a_1},\dots,K_{a_k}$,\, $a_0+a_1+\dots+a_k=n$. Without loss of generality, assume that $a_0\geq a_1\geq a_2\geq\dots\geq a_k\geq1$.

If $a_0=1$, then $G$ is an $n$-vertex tree. By [13, Theorems 3.3 and 3.8], we know the second last tree in an $S$-order among $n$-vertex trees is just $K_{1,n-1}^*$. It is easy to see that $a_0\not=2,$ hence in what follows we consider $a_0\ge 3.$

Let $V_i=\{v\in V_{K_{a_i}}\!:\, v$ is an end-vertex of a cut edge of $G$\}. In order to complete the proof, it suffices to show the following facts.

\vspace{2mm}\noindent {\bf Fact 1}. If $a_0\geq3$, then $|V_0|=2, |V_1|=|V_2|=\cdots=|V_k|=1$.
\begin{proof}
We prove Fact 1 by contradiction. If $|V_0|=|V_1|=|V_2|=\cdots=|V_k|=1$, then without loss of generality assume that $V_i=\{u_i\}$,\, $i=0,1,\ldots, k$.

First we consider $G\in \mathscr{K}_n^k$. Note that $G\in \mathscr{K}_n^k\setminus\{K_n^k\}$, hence $a_1\geq3$; otherwise, $a_1=2$, which implies that $G$ contains at least $k+1$ cut edges, a contradiction. If $a_1>3$, we consider graph $G^*:=K(a_0+1,\{a_1-1, a_2,\ldots, a_k\})$ in $\mathscr{K}_n^k\setminus\{K_n^k\}$.
Note that $S_i(G)=S_i(G^*)$ for $i=0,1$ and $S_2(G)-S_2(G^*)=2(a_1-1-a_0)<0$, hence $G\prec_sG^*$, a contradiction. Therefore $a_1=3$.

If $a_2>1$, we consider graph $G':=K(a_0+a_2-1,\{3, 1,a_3,\ldots, a_k\})\in \mathscr{G}_n^k\backslash K_n^k$. Note that $ S_i(G)=S_i(G')$ for $i=0,1$ and $S_2(G)-S_2(G')=2(a_2-1-(a_2-1)a_0)=2(a_2-1)(1-a_0)<0$, hence $G\prec_s G'$, a contradiction. Therefore, $a_2=1$, whence $a_3=\dots=a_k=1$. Together with $a_1=3$, we have $a_0=n-k-2.$ That is to say, $G\cong K(n-k-2,\{3,1,1,\ldots, 1\}).$

For convenience, let $w_1\in N_{K_{a_0}}(u_0)$ and $N_G(u_1)=\{u_0,v_1,v_2\}.$ Consider
$$
G^*:=G-\{u_1v_1,u_1v_2\}+\{w_1v_1,w_1v_2\},
$$
it is easy to see that $G^*\in\mathscr{G}_n^k\backslash K_n^k$. Note that $S_i(G)=S_i(G^*)$ for $i=0,1,2,3$ and $\phi_G(P_2)=\phi_{G^*}(P_2), \phi_G(C_4)=\phi_{G^*}(C_4)$,
$\phi_G(P_3)-\phi_{G^*}(P_3)=2-2(a_0-1)=2(2-a_0)<0$, hence by Lemma 1.2(i) $S_4(G)<S_4(G^*)$. Thus $G\prec_s G^*$, a contradiction.
Therefore $G\not\in\mathscr{K}_n^k$.

Now we consider the case $G\not\in\mathscr{K}_n^k$. It is easy to see that the edge induced graph $G[\Bbb{E}]$ is a tree which is not isomorphic to $K_{1,k}.$
Hence, partition $V_{G[\Bbb{E}]}$ into $D^0(G[\Bbb{E}])\cup D^1(G[\Bbb{E}])\cup D^2(G[\Bbb{E}])\cup D^3(G[\Bbb{E}])\cup \cdots$, where $D^i(G[\Bbb{E}])=\{u\in V_{G[\Bbb{E}]}: \dist_{G[\Bbb{E}]}(u, u_0)=i\},\, i=0,1,2,3,\ldots.$
It is easy to see that $D^2(G[\Bbb{E}])\not=\emptyset.$

If $D^3(G[\Bbb{E}])\not=\emptyset$, that is to say, there exists $u\in D^2(G[\Bbb{E}])$ such that $d_{G[\Bbb{E}]}(u)\ge 2,$ then choose $u_i$ from $D^1(G[\Bbb{E}])$ such that $u_i$ is adjacent to $u_0$ and $u$. Let $W:=N_{G[\Bbb{E}]}(u_i)\setminus \{u_0\}.$ As $u\in W$, we have $W\not=\emptyset$. Consider
$$
G^*=G-\{u_iw:\, w\in W\}+\{u_0w:\, w\in W\},
$$
then its routine to check that $G^*\in\mathscr{G}_n^k\setminus \{K_n^k\}$. Note that $S_i(G)=S_i(G^*)$ for $i=0, 1, 2,3$ and $\phi_G(P_2)=\phi_{G^*}(P_2), \phi_G(C_4)=\phi_{G^*}(C_4)$, hence by Lemma 1.2(i) we have
$$
  S_4(G)-S_4(G^*) = 4(\phi_G(P_3)-\phi_{G^*}(P_3))=4[(a_i-a_0)-st],
$$
where $s=|W|\ge 1$ and $t=|N_{G[\Bbb{E}]}(u_0)\setminus\{u_i\}|\ge 0$.
Note that $a_i\le a_0$, hence if $a_i<a_0$, then for all $t\ge 0$ we have $(a_i-a_0)-st<0$, which implies that $G\prec_sG^*$, a contradiction. If
$a_i=a_0$, then for all $t\ge 1$ we have $(a_i-a_0)-st<0$, which implies that $G\prec_sG^*$, a contradiction. If $a_i=a_0$ and $t = 0$, then
$G^*\cong G$. Hence, in order to complete the proof, it suffices to consider $D^3(G[\Bbb{E}])=\emptyset$ and $d_{G[\Bbb{E}]}(u_0)>1.$ Furthermore, as $G\not\in \mathscr{K}_n^k$, we have $D^2(G[\Bbb{E}])\not=\emptyset$ and for each $u\in D^2(G[\Bbb{E}])$, $u$ is a leaf of $G[\Bbb{E}]$ (otherwise, $D^3(G[\Bbb{E}])\not=\emptyset,$ a contradiction).

If there exists $u_i\in D^1(G[\Bbb{E}])$ such that $d_{G[\Bbb{E}]}(u_i)\ge 3$, then move $d_{G[\Bbb{E}]}(u_i)-2$ pendant edges to $u_0$ and denote the resultant graph by
$G'$. It is easy to see that $G'\in\mathscr{G}_n^k \setminus \{K_n^k\}.$ Note that $s:=d_{G[\Bbb{E}]}(u_i)-2\ge 1,\, q:=d_{G[\Bbb{E}]}(u_0)-1\ge 1,\, S_i(G)=S_i(G')$ for $ i=0,1,2,3$ and $\phi_G(P_2)=\phi_{G'}(P_2),\phi_G(C_4)=\phi_{G'}(C_4)$, hence by Lemma~1.2(i) we have
\begin{eqnarray}
S_4(G)-S_4(G')    &=& 4(\phi_G(P_3)-\phi_{G'}(P_3))\notag\\
                  &=&  4[s(a_i-1)-s(a_0-1)-s(q-1)]\notag\\
                  &=& 4s(a_i-a_0-q+1).
\end{eqnarray}
If $a_0>a_i$ or $q\geq2$, in the view of (2.1), we obtain that $S_4(G)-S_4(G')<0$, i.e., $G\prec_sG'$, a contradiction. If  $a_0=a_i$ and $q=1$, then it is easy to see $G'\cong G$. Hence, in order to complete the proof, it suffices to consider that, in the edge induced graph $G[\Bbb{E}]$, each of the non-pendant vertices in $D^1(G[\Bbb{E}])$ is of degree 2.

For convenience, let $W=\{u: u\in D^1(G[\Bbb{E}]), d_{G[\Bbb{E}]}(u)=2\}.$ It is easy to see that $W\not=\emptyset.$ If $|W|\ge 2$, choose $u\in W$ such that
its unique neighbor in $G[\Bbb{E}]$ is a leaf, say $u'$. Let $$
G^*=G-uu'+u_0u',
$$
then $G^*\in\mathscr{G}_n^k\setminus\{K_n^k\}$.
Since $S_i(G)=S_i(G^*)$ for $ i=0,1,2,3$ and $\phi_G(P_2)=\phi_{G^*}(P_2),\phi_G(C_4)=\phi_{G^*}(C_4)$,
$$
\phi_G(P_3)-\phi_{G^*}(P_3)=(a_i-1)-(a_0-1)-p=a_i-a_0-p<0,
$$
we have $S_4(G)-S_4(G^*)<0,$ i.e., $G\prec_sG^*$, a contradiction.
Hence, $|W|=1$.

By a similar discussion as in the proof of Fact 3 in Theorem 2.1, we can obtain that $a_0=n-k, a_1=a_2=\ldots=a_k=1$. Note that $a_0\geq 3$, hence $k<n-1$. Assume that $W=\{u\}$ with $N_{G[\Bbb{E}]}(u)=\{u_0, u'\}$, where $u'$ is a pendant vertex in $G[\Bbb{E}]$.
Let $x\in N_{K_{a_0}}(u_0)$. Consider
$
G^*=G-\{uu'\}+\{xu'\},
$
then $G^*\in\mathscr{G}_n^k\setminus\{K_n^k\}$.
Note that $S_i(G)=S_i(G^*)$ for $i=0,1,2,3$,\, $\phi_G(P_2)=\phi_{G^*}(P_2), \phi_G(C_4)=\phi_{G^*}(C_4)$ and $\phi_G(P_3)-\phi_{G^*}(P_3)=1-(a_0-1)=2-a_0<0$ $(a_0\geq3)$, hence $S_4(G)-S_4(G^*)<0,$ i.e., $G\prec_sG^*$, a contradiction.

If there is a $V_i$ satisfying $|V_i|\geq3$, then choose two distinct vertices $u_i', u_i''$ in $V_i$ and let
$
G^*=G-\{u_i'u:\, u\in N_{G[\Bbb{E}]}(u_i')\}+\{u_i''u:\, u\in N_{G[\Bbb{E}]}(u_i')\}.
$
It is easy to see that $G^*\in\mathscr{G}_n^k\setminus\{K_n^k\}$. Note that $S_i(G)=S_i(G^*), i=0, 1, 2, 3$, \, $\phi_G(P_2)=\phi_{G^*}(P_2),\phi_G(C_4)=\phi_{G^*}(C_4)$ and $\phi_G(P_3)-\phi_{G^*}(P_3)=-|N_{G[\Bbb{E}]}(u_i')||N_{G[\Bbb{E}]}(u_i'')|<0$, hence $S_4(G)-S_4(G^*)<0,$ which implies that $G\prec_sG^*$, a contradiction.\vspace{2mm}

If there exists an $i\in \{1,2,\ldots, k\}$ such that $|V_i|=2$. Assume, without loss of generality, that
$V_i=\{u_i, u'\}$, where $u_i$ is adjacent to $u_0\in V_0$. Let
$$
G^*=G-\{u_ix: x\in V_{K_{a_i}}\}+\{yx: y\in V_{K_{a_0}}, x\in V_{K_{a_i}}\setminus\{u_i\}\}.
$$
It is easy to see that $G^*\in\mathscr{G}_n^k\backslash K_n^k$. Note that $S_i(G)=S_i(G^*)$ for $i=0,1$ and $S_2(G)-S_2(G^*)=2(a_i-1)-2(a_i-1)a_0=2(a_i-1)(1-a_0)<0$, hence $S_2(G)<S_2(G^*),$ i.e., $G\prec_s G^*$, a contradiction.

Combining with discussion as above, we obtain that $|V_1|=|V_2|=\cdots=|V_k|=1,$ whence $|V_0|=2$, as desired.
\end{proof}


\vspace{2mm}\noindent {\bf Fact 2}. $a_1=a_2=\cdots=a_k=1$.
\begin{proof}
By a similar discussion as in the proof of Fact 3 in Theorem 2.1, we can get $a_0=n-k$, $a_1=a_2=\dots=a_k=1$. We omit the procedure here.
\end{proof}

\vspace{2mm}\noindent {\bf Fact 3}. $G\cong G_3$, where $G_3=K(n-k,\{1,1,\dots,1\})-u_0u_k+u_0'u_k$, where $u_0u_k$ is a cut edge and $u_0'\in V_{K_{a_0}}\backslash \{u_0\}$.
\begin{proof}
Note that if $G$ has just two cut edges, it is easy to see that $G\cong G_3$ defined as above. Hence in what follows we consider that $G$ contains at least three cut edges.

Let $N_{G[\Bbb{E}]}(u_0)=\{u_1,u_2,\dots,u_m\}$ and $N_{G[\Bbb{E}]}(u_0')=\{u_1',u_2',\dots,u_t'\}$. Without loss of generality, assume that $m\geq t$. Obviously, $t\ge 1$.
At first we show that $t=1$. Otherwise, let
$$
     G^*=G-\{u_0'u_2',u_0'u_3',\dots,u_0'u_t'\}+\{u_0u_2',u_0u_3',\dots,u_0u_t'\}.
$$
It is easy to see that $G^*\in\mathscr{G}_n^k\backslash K_n^k$.
Note that $S_i(G)=S_i(G^*)$ for $ i=0,1,2,3$,\, $\phi_G(P_2)=\phi_{G^*}(P_2),$ and $\phi_G(C_4)=\phi_{G^*}(C_4)$, hence
$$
S_4(G)-S_4(G^*)=4(\phi_G(P_3)-\phi_{G^*}(P_3))=-4m(t-1)<0,
$$
i.e., $G\prec_sG^*$, a contradiction.

Now we are to show that $m=k-1.$ If not, there exists a vertex $u\in \{u_1,u_2,\dots,u_m, u_1',u_2',\dots,u_t'\}$ such that $d_{G[\Bbb{E}]}(u)\ge 2.$ Denote $N_{G[\Bbb{E}]}(u)\setminus\{u_0,v_0\}=\{\hat{u}_1,\hat{u}_2,\dots,\hat{u}_s\}$, $s\geq1$. Let
$$
G^*=G-\{u\hat{u}_1,u\hat{u}_2,\dots,u\hat{u}_s\}+\{u_0\hat{u}_1,u_0\hat{u}_2,\dots,u_0\hat{u}_s\}.
$$
It is easy to see that $G^*\in\mathscr{G}_n^k\backslash K_n^k$.
Notice that $S_i(G)=S_i(G^*)$ for $ i=0,1,2,3$,  $\phi_G(P_2)=\phi_{G^*}(P_2)$ and $\phi_G(C_4)=\phi_{G^*}(C_4)$, hence
$$
S_4(G)-S_4(G^*)=4(\phi_G(P_3)-\phi_{G^*}(P_3))=4s((a_i-1)-(a_0-1)-(m-1))=4s(a_i-a_0)-4s(m-1)<0.
$$
The last inequality follows by $a_i=1< n-k = a_0$ (by fact 2), and $m\ge 1, s\ge 1$. Hence, we get $S_4(G)<S_4(G^*),$ i.e., $G\prec_sG^*$, a contradiction.
so we have $m=k-1, t=1$, which is equivalent to that $G\cong G_3$.
\end{proof}

This completes the proof.
\end{proof}

\section{\normalsize The first and the second graphs in an $S$-order among $\mathscr{G}_n^k$ }\setcounter{equation}{0}
In this section, we are to determine the first and the second graphs in an $S$-order among $\mathscr{G}_n^k$.
Let ${\Bbb E}=\{e_1,e_2,\dots,e_k\}$ be the set of all the cut edges of $G\in\mathscr{G}_n^k$. Note that if we delete an edge, say $e$, from a connected graph $G$, then in the view of $S_2(G)=2|E_G|$, we have $S_2(G)>S_2(G-e)$. In order to determine the first graph in an $S$-order among $\mathscr{G}_n^k$, it suffices to choose the graph such that its size is as small as possible.
\begin{thm}
Of all the connected graphs with $n$ vertices and $k$ cut edges, the first graph in an $S$-order is obtained uniquely at $P_n^k$.
\end{thm}
\begin{proof}
Choose $G\in \mathscr{G}_n^k$ such that it is as small as possible according to the relation $\preceq_s$. If $k=0$, then it is easy to see that $G\cong C_n$ and our result holds immediately. Therefore we may assume that $k\geq1$. We show the following claim at first.
\begin{claim}
$G$ contains exactly one cycle.
\end{claim}
\begin{proof}
Assume to the contrary that $G$ contains at least two cycles. If $G$ contains two cycles $C^1$ and $C^2$ such that $C^1$ and $C^2$ have edges in common; see Fig. 2(a), then let $G^*=G-\{uv, xy\}+ux$ (see Fig. 2(b)); if $G$ contains two cycles $C^1$ and $C^2$ such that $C^1$ and $C^2$ have just one vertex in common; see Fig. 2(c), then let $G^*=G-\{ux, vx\}+uv$ (see Fig. 2(d)). It is routine to check that $G^*\in \mathscr{G}_n^k$
\begin{figure}[h!]
\begin{center}
\psfrag{a}{$C^1$} \psfrag{b}{$C^2$}
\psfrag{t}{$C'$} \psfrag{r}{$C''$}
\psfrag{c}{$u$} \psfrag{d}{$v$}
\psfrag{e}{$x$} \psfrag{f}{$y$}
\psfrag{1}{(a)} \psfrag{2}{(b)}
\psfrag{3}{(c)} \psfrag{4}{(d)}
 \includegraphics[width=120mm]{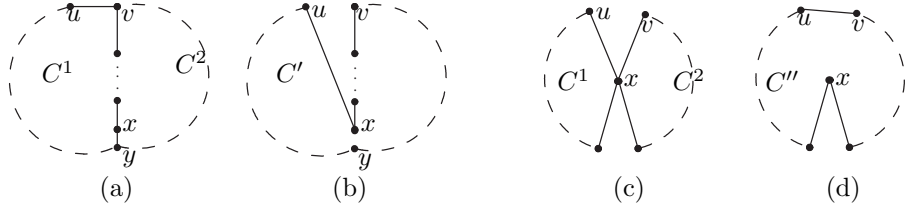}\\
  \caption{Graphs used in the proof of Claim 1.}
\end{center}
\end{figure}
and in each of the above cases one has $S_i(G)=S_i(G^*),\,i=0,1$ and $S_2(G)-S_2(G^*)=2>0$, hence $ G^*\prec_s G$, a contradiction.

If $G$ contains two cycles $C_l=u_0u_1u_2\dots u_{l-1} $ and $C_j=v_0v_1v_2\dots v_{j-1}$ such that $C_l$ connects $C_j$ by a path $P_i,\, i\ge 2$, whose end vertices are $u_0, v_1$, and the vertex, say $u_t$ (resp. $v_m$), on the cycle  $C_l$ (resp. $C_j$) in $G$ either is of degree 2 or has subgraph $G_t$ (resp. $H_m$) attached, $ 0\leq t \leq {l-1} $, $ 0\leq m \leq {j-1} $; see Fig. 3. Let
$$
  G^*=G-\{u_0u_1,v_1v_2,v_0v_1\}+\{u_0v_2,u_1v_0\} ,
$$
then $G^*\in\mathscr{G}_n^k$. Since $S_i(G)=S_i(G^*),\,i=0,1$. $S_2(G)-S_2(G^*)=2>0$, then $ G^*\prec_s G$, a contradiction. Therefore, $G$ contains exactly one cycle.
\end{proof}
\begin{figure}[h!]
\begin{center}
 \includegraphics[width=150mm]{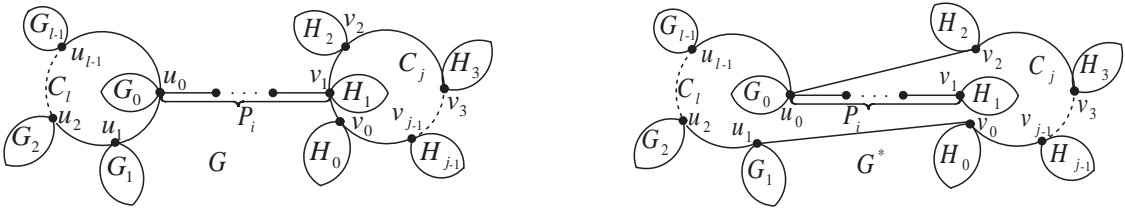}\\
  \caption{Graph  $G\Rightarrow G^*$.}
\end{center}
\end{figure}

By Claim 1, we know that $G$ is a unicyclic graph. Note that $G$ contains exactly $k$ cut edges, hence $G$ is an $n$-vertex unicyclic graph with girth $n-k$.
By Lemma \ref{lem4} the first graph in an $S$-order among the $n$-vertex unicyclic graph with girth $n-k$ is just the graph $P_n^k,$ as desired.
\end{proof}

At the rest of this section, we are to determine the second graph in an $S$-order among $\mathscr{G}_n^k\, (k \geq 3)$.
\begin{thm}
Of all graphs with $n$ vertices and $k$ cut edges, the second graph in an $S$-order is obtained uniquely at $\hat{U}_n^k\, (k\geq 3)$, where $\hat{U}_n^k$ is obtained by attaching two leafs to the pendant vertex of graph $P_{n-2}^{k-2}$.
\end{thm}
\begin{proof}
Note that if we delete an edge $e$ from a connected graph $G$, then in the view of $S_2(G)=2|E_G|$, we have $S_2(G)>S_2(G-e)$, hence in order to determine the second graph in an $S$-order among $\mathscr{G}_n^k$, it suffices to determine the second graph in an $S$-order among the set of all $n$-vertex unicyclic graphs with girth $n-k$; we denote this set by $\mathscr{U}_n^k$.

Choose $G\in \mathscr{U}_n^k\setminus\{P_n^k\}$ such that it is as small as possible with respect to $\preceq_s$. Note that $\Bbb{E}$ is the set of $k$ cut edges of $G$, hence $G[\Bbb{E}]$ is a forest. We are to show that $G[\Bbb{E}]$ is a tree. If this is not true, then it is equivalent to that there exist at least two vertices, say $u_0, v_0$, on the unique cycle contained in $G$ satisfying $d_G(u_0), d_G(v_0)\ge 3.$

In the edge induced graph $G[\Bbb{E}]$, consider the tree, say $T_1$, containing $u_0$. We are to show that $T_1$ is a path; otherwise, choose a longest path $P=u_0u_1\dots u_p$ in $T_1$ with end-vertex $u_0, u_p$, it is easy to see $d_{T_1}(u_p)=1$.
If there exists $u_i$ with $i\ge 1$ on $P$ such that $d_G(u_i)>2$. Choose a vertex $x$ in $N_G(u_i)\setminus\{u_{i-1}, u_{i+1}\}$ and let
$
G^*=G-u_ix+u_px,
$
then $G^*\in\mathscr{G}_n^k\backslash \{P_n^k\}$.
Note that $S_i(G)=S_i(G^*)$ for $ i=0,1,2,3$,\, $\phi_G(P_2)=\phi_{G^*}(P_2),\phi_G(C_4)=\phi_{G^*}(C_4)$ and $\phi_G(P_3)-\phi_{G^*}(P_3)\geq1$, hence by
Lemma~1.2(i), we get $S_4(G)-S_4(G^*)>0$, i.e., $G^*\prec_s G$, a contradiction. Hence, we obtain that each vertex $u_i$ on $P$ is of degree 2 in $G$ for  $i=1,2,\ldots, p-1$. Hence, if $d_G(u_0)=3$, then $T_1$ is a path, as desired.
If $d_G(u_0)>3$, then choose $x$ from $N_G(u_0)$ such that $x$ is not on the cycle and the path $P$ contained in $G$. Let
$
G^*=G-u_0x+u_px,
$
then $G^*\in\mathscr{G}_n^k\backslash P_n^k$.
Notice that $S_i(G)=S_i(G^*)$ for $ i=0,1,2,3$, $\phi_G(P_2)=\phi_{G^*}(P_2),\phi_G(C_4)=\phi_{G^*}(C_4)$ and $\phi_G(P_3)-\phi_{G^*}(P_3)\geq\emph{}2$,  by Lemma 1.2(i), we get $S_4(G)-S_4(G^*)>0,$ i,e., $G^*\prec_s G$, a contradiction.
By a similar discussion as above, we can also show that, in $G[\Bbb{E}]$, the component contains $v_0$ is also a path, say $P'$. For convenience, let $v_0'$ be the neighbor of $v_0$ on $P'$.

   If there exists another vertex $u_0'\neq u_0, v_0$, on the unique cycle contained in $G$ satisfying $d_G(u_0')\ge 3$. Let $G^*=G-\{u_0'x, x\in N_{G[\Bbb{E}]}(u_0')\}+\{u_px, x\in N_{G[\Bbb{E}]}(u_0')\}$,
then $G^*\in\mathscr{G}_n^k\backslash P_n^k$. Note that $S_i(G)=S_i(G^*)$ for $ i=0,1,2,3$, $\phi_G(P_2)=\phi_{G^*}(P_2),\phi_G(C_4)=\phi_{G^*}(C_4)$ and $\phi_G(P_3)-\phi_{G^*}(P_3)\geq\emph{}1$,  by Lemma 1.2(i), we get $S_4(G)-S_4(G^*)>0,$ i,e., $G^*\prec_s G$, a contradiction. So, we just need to consider there exist two vertices
on the unique cycle contained in $G$ satisfying $d_G(u_0), d_G(v_0)\ge 3.$ Without loss of generality, assume that $|E_P|\geq |E_{P'}|$. Let $G^*=G-v_0v_0'+u_{p-1}v_0'$, it is easy to see that
$G^*\in\mathscr{G}_n^k\setminus\{P_n^k\}$.

\noindent $\bullet$ $k=3$. By Lemma 1.1, we have $S_i(G)\geq S_i(G^*)$ for $i=0,1,\dots,n-2$ and
$S_{n-1}(G)>S_{n-1}(G^*)$. Hence $G^*\prec_sG$, a contradiction.

\noindent $\bullet$ $k\geq4$. Note that $\phi_G(P_2)=\phi_{G^*}(P_2),\phi_G(C_4)=\phi_{G^*}(C_4)$, $\phi_G(P_3)=\phi_{G^*}(P_3),\phi_G(K_{1,3})=\phi_{G^*}(K_{1,3}),\linebreak \phi_G(U_5)=\phi_{G^*}(U_5), \phi_G(B_4)=\phi_{G^*}(B_4)$,
$\phi_G(B_5)=\phi_{G^*}(B_5)$, $\phi_G(C_3)=\phi_{G^*}(C_3)$, $\phi_G(C_6)=\phi_{G^*}(C_6)
$ and $\phi_G(P_4)-\phi_{G^*}(P_4)\geq1$, hence by Lemma 1.2, 1.3, we get that $S_i(G)=S_i(G^*)$ for $i=0,1,2,3,4,5$ and $S_6(G)-S_6(G^*)>0,$ i.e., $G^*\prec_sG$, a contradiction.

Therefore, we obtain that $G[\Bbb{E}]$ is a tree. That is to say, there exists just one vertex, say $u_0$, on the unique cycle such that $d_G(u_0)\ge 3$.
Choose one of the longest paths, say $P:=u_0u_1\dots u_p$, from $G[\Bbb{E}]$. It is easy to see that $u_p$ is a leaf of $G$. Furthermore, we have the following claim.
\begin{claim}
The length of $P$ is $k-1$, i.e., $P:=u_0u_1\dots u_{k-2}u_{k-1}$ and $G[\Bbb{E}]$ is obtained from $P$ by attaching a leaf to $u_{k-2}$ of $P$.
\end{claim}
\begin{proof}
Note that $P=u_0u_1\dots u_p$ is one of the longest paths of $G[\Bbb{E}]$ and $u_p$ is a leaf. Hence, we first show that $d_G(u_0)=3$. Otherwise, choose $x$ from $N_G(u_0)$ such that $x$ is not on the cycle and the path $P$ of $G$. If $d_G(u_i)\geq3$ for some $i\in \{1,2,\dots,p-1\}$, let
$
G^*=G-u_0x+u_px.
$
Obviously, $G^*\in\mathscr{G}_n^k\backslash\{P_n^k\}$.
Note that $S_i(G)=S_i(G^*)$, $ i=0,1,2,3$,\, $\phi_G(P_2)=\phi_{G^*}(P_2),\phi_G(C_4)=\phi_{G^*}(C_4)$ and $\phi_G(P_3)-\phi_{G^*}(P_3)\geq 2$, hence by Lemma~1.2(i), we get $S_4(G)-S_4(G^*)>0,$ i.e., $G^*\prec_s G$, a contradiction. If $d_G(u_i)=2$ for any $i\in \{1,2,\dots,p-1\}$, let
$
G^*=G-u_0x+u_{p-1}x.
$
Obviously, $G^*\in\mathscr{G}_n^k\backslash\{P_n^k\}$.
Note that $S_i(G)=S_i(G^*)$, $ i=0,1,2,3$,\, $\phi_G(P_2)=\phi_{G^*}(P_2),\phi_G(C_4)=\phi_{G^*}(C_4)$ and $\phi_G(P_3)-\phi_{G^*}(P_3)\geq 1$, hence by Lemma~1.2(i), we get $S_4(G)-S_4(G^*)>0,$ i.e., $G^*\prec_s G$, a contradiction.

Now we show that $d_G(u_i)=2$, $i=1,2,\dots,p-2$ and $d_G(u_{p-1})=3$.
Note that $G\ncong P_n^k$, hence there exists at least one vertex $u_i\, (1 \leq i\leq {p-1})$ on $P$ such that $d_G(u_i)\geq3$.

If there exists a vertex $u_i\, (1 \leq i\leq {p-1})$ on $P$ such that $d_G(u_i)\geq 4$, then choose $x\in N_G(u_i)\setminus\{u_{i-1}, u_{i+1}\}$ and let
$
G^*=G-u_ix+u_px.
$
Obviously, $G^*\in\mathscr{G}_n^k\backslash P_n^k$. Since $S_i(G)=S_i(G^*)$, $ i=0,1,2,3$. $\phi_G(P_2)=\phi_{G^*}(P_2),\,\phi_G(C_4)=\phi_{G^*}(C_4)$ and $\phi_G(P_3)-\phi_{G^*}(P_3)\geq2$, by Lemma 1.2(i), we have $S_4(G)-S_4(G^*)>0,$ i.e., $G^*\prec_sG$, a contradiction. Hence,  $\max\{d_G(u_i),i=1,2,\dots,p-1\}=3$.


If $d_G(u_{p-1})=d_G(u_i)=3$ for some $i\in \{1, 2, \ldots, p-2\}$, then choose  $z_1$ in $N_G(u_i)\setminus\{u_{i-1},u_{i+1}\}$ and let
$
G^*=G-u_iz_1+u_pz_1,
$
then $G^*\in\mathscr{G}_n^k\backslash P_n^k$.
Notice that $S_i(G)=S_i(G^*)$, $i=0,1,2,3$. $\phi_G(P_2)=\phi_{G^*}(P_2),\,\phi_G(C_4)=\phi_{G^*}(C_4)$ and $\phi_G(P_3)-\phi_{G^*}(P_3)=1>0$, hence by Lemma 1.2(i), we get that $S_4(G)-S_4(G^*)>0,$ i.e., $G^*\prec_sG$, a contradiction.

If $d_G(u_{p-1})=2,\, d_G(u_i)=3$ for some $i\in \{1, 2, \ldots, p-2\}$, then choose  $z_1$ in $N_G(u_i)\setminus\{u_{i-1},u_{i+1}\}$ and let
$
G^*=G-u_iz_1+u_{p-1}z_1,
$
then it is easy to see that $G^*\in\mathscr{G}_n^k\backslash \{P_n^k\}$.
Note that $S_i(G)=S_i(G^*)$ for $i=0,1,2,3,4,5$, $\phi_G(P_2)=\phi_{G^*}(P_2),\,\phi_G(C_4)=\phi_{G^*}(C_4)$, $\phi_G(P_3)=\phi_{G^*}(P_3),\,\phi_G(K_{1,3})=\phi_{G^*}(K_{1,3}),\,\phi_G(U_5)=\phi_{G^*}(U_5),\, \phi_G(B_4)=\phi_{G^*}(B_4),$ $\phi_G(B_5)=\phi_{G^*}(B_5),\,\phi_G(C_3)=\phi_{G^*}(C_3),\,$
$\phi_G(C_6)=\phi_{G^*}(C_6)$ and $\phi_G(P_4)-\phi_{G^*}(P_4)\geq1$, hence by Lemma 1.2(iii), we get that $S_6(G)-S_6(G^*)>0,$ i.e., $G^*\prec_sG$, a contradiction.

Hence, we obtain that $d_G(u_0)=3, d_G(u_1)=d_G(u_2)=\ldots=d_G(u_{p-2})=2, d_G(u_{p-1})=3$ and $d_G(u_p)=1.$ For convenience, let $N_G(u_{p-1})\setminus\{u_{p-2},u_p\}=\{z_0\}$. It is easy to see that $z_0$ is a leaf; otherwise $G[\Bbb{E}]$ contains a path $P':=u_0u_1\ldots u_{p-1}z_0\ldots z_t$, where $z_t$ is a leaf. It is routine to check that the length of $P'$ is longer than that of $P$, a contradiction. Therefore, $d_G(z_0)=1$, as desired.
\end{proof}
Based on Claim 2, Theorem 3.2 follows immediately.
\end{proof}

\end{document}